# On a new class of fractional partial differential equations

Tien-Tsan Shieh and Daniel E. Spector

*Dedicated to Israel and Taiwan*

**Abstract.** In this paper, we study a new class of fractional partial differential equations which are obtained by minimizing variational problems in fractional Sobolev spaces. We introduce a notion of fractional gradient which has the potential to extend to many classical results in the Sobolev spaces to the nonlocal and fractional setting in a natural way.

## 1 Introduction and Main Results

Introduced by Marcel Riesz [25–27] in the context of potential theory, the fractional Laplacian is ubiquitous in the modern study of fractional partial differential equations. It enjoys a variety of definitions, though they can be distilled down to the following two: as a pseudodifferential operator in Fourier space[1]

$$(-\Delta)^s u := ((2\pi|\xi|)^{2s}\widehat{u})^{\vee} \qquad (1)$$

or as a singular integral in real space

$$(-\Delta)^s u(x) := c_{N,s}\ p.v. \int_{\mathbb{R}^N} \frac{u(x)-u(y)}{|x-y|^{N+2s}}\ dy. \qquad (2)$$

Here, and in the sequel, we restrict ourselves to the consideration of $s \in (0,1)$, though in general, definitions involving spectral representations fall in the former camp, while those involving derivation of potentials fall in the latter. We recall

The first author is supported in part by National Science Council of Taiwan under research grant NSC-101-2115-M-009-014-MY3. The second author is supported in part by a Technion Fellowship.

[1] Here we use the notation $\widehat{u}$ for the Fourier transform and $\widecheck{u}$ for its inverse, with the convention that $\widehat{u}(\xi) = \int_{\mathbb{R}^N} u(x)e^{-2\pi i x \cdot \xi}\ dx$.

that Riesz had introduced generalized potentials of order $\alpha$ (now called Riesz potentials) for $0 < \alpha < N$ by the formula[2]

$$I_\alpha u := I_\alpha * u,$$

where

$$I_\alpha(x) := \frac{\gamma(N,\alpha)}{|x|^{N-\alpha}},$$

and the constant

$$\gamma(N,\alpha) := \frac{\Gamma(\frac{N-\alpha}{2})}{\pi^{\frac{N}{2}} 2^\alpha \Gamma(\frac{\alpha}{2})}.$$

Regarding these potentials, among his fundamental results are that they satisfy the semi-group property,

$$I_\alpha I_\beta u = I_{\alpha+\beta} u, \tag{3}$$

for $\alpha, \beta > 0$ and $\alpha + \beta < N$, and that the Laplacian maps a potential of order $\alpha + 2$ to a potential of order $\alpha$,

$$-\Delta I_{\alpha+2} u = I_\alpha u. \tag{4}$$

Through analytic continuation, the Riesz potential can be extended to negative exponents (see Chapter I, p.41 in [21] for a detailed treatment), and thus one arrives at a formula for the fractional Laplacian, or Riesz fractional derivative, that is in agreement with (2),

$$(-\Delta)^s u := I_{-2s} u. \tag{5}$$

Classical interest in the fractional Laplacian had been through its definition as a pseudodifferential operator (see, for instance, the work of Bochner [5], Yoshida [34], and Kato [19]), and while its spectral representation remains an important tool in understanding questions of regularity, it is often taken in tandem with geometric-space estimates based on the singular integral representation. For example, in recent work of Da Lio and Rivière [8], Da Lio and Schikorra [9], and Schikorra [28–30] we find alternative use of the two in establishing Hölder continuity for critical points of integral functionals of the fractional Laplacian where the underlying fields take values in a manifold. From the standpoint of fractional partial differential equations, there has been a resurgence of interest in its definition as a Cauchy principle value integral, since this gives a prototype for general singular integral operators (or nonlocal operators, as they are termed by the community working on them). We can cite, for instance, the regularity theory of nonlocal equations developed by Fukushima [13] and Kassmann [17, 18] in the symmetric case, and Barles, Chasseigne, and Imbert [3], Caffarelli and Silvestre [7], and Di Castro, Kuusi, and Palatucci [11] where no such symmetry is assumed.

---

[2] Riesz had introduced generalized potentials by the formula $I^\alpha$, though we prefer to follow the notation in Chapter V, p.117 of Stein [31] using $I_\alpha$.

That the fractional Laplacian is the central object in such a diverse set of problems is interesting when contrasted with the non-fractional setting, since though the (negative) Laplacian is the prototype of an elliptic operator, it is not the foundational object of the calculus of variations or partial differential equations. It is this observation that is the starting point of this paper, along with the assertion that there is an object more fundamental than the fractional Laplacian. This object, the distributional Riesz fractional gradient, will shortly be defined, and we will show that it is the intrinsic object of interest for the study of fractional problems in the calculus of variations and partial differential equations.

**Definition 1.1.** Let $s \in (0, 1)$. If $u \in L^p(\mathbb{R}^N)$ for some $1 < p < \infty$ such that $I_{1-s} * u$ is well-defined, we define the *distributional Riesz fractional gradient*

$$(D^s u)_j := \frac{\partial^s u}{\partial x_j^s}, \quad j = 1, \ldots, N,$$

where

$$\frac{\partial^s u}{\partial x_j^s} := \frac{\partial}{\partial x_j} I_{1-s} u, \tag{6}$$

in the sense that

$$\langle \frac{\partial^s u}{\partial x_j^s}, v \rangle = (-1) \langle I_{1-s} u, \frac{\partial v}{\partial x_j} \rangle,$$
$$= -\int_{\mathbb{R}^N} (I_{1-s} * u) \frac{\partial v}{\partial x_j} \, dx$$

for every $v \in C_c^\infty(\mathbb{R}^N)$.

The importance of such an object will soon become evident, though it can be summarized in its properties analogous to the gradient: it has a geometric interpretation, behaves well under the distributional calculus, and gives rise to spaces of weakly differentiable functions. Thus, although the germ of such an idea is present in the work of Riesz [27], the results we detail would not be possible without the subsequent development of Sobolev spaces, potential theory, and harmonic analysis. In this respect, it would be surprising if such an important object had escaped all interest of researchers. In fact, it has been mentioned without further study in some existing literature - an equivalent definition[3] of the distributional Riesz fractional gradient has been given by Schikorra [29], while it can also be derived as a special case of the nonlocal gradients of Gilboa and Osher [16]. In this paper, we will give a full exposition of the distributional Riesz fractional gradient, detailing its properties, as well as to pose a new class of fractional problems in the calculus of variations and partial differential equations.

We have chosen the name distributional Riesz fractional gradient to acknowledge Riesz's contribution through his definition of generalized potentials, as well as to distinguish it from other fractional derivatives. Indeed, the literature on fractional derivatives in one dimension is quite extensive[4], while in several dimensions there have been other possible definitions due to Meerschaert, Mortensen, and Wheatcraft [23] or Tarasov [33]. What distinguishes the distributional Riesz fractional gradient from existing notions is the rich theory we develop for spaces of functions possessing a fractional derivative, made possible by its linearity and reliance on the distributional calculus. Connections of this nature have been made for **nonlinear** fractional differential objects, with early contributions due to Stein and Strichartz in the context of Bessel potential spaces (see [31] pp.161-163 for extensive references), while more recent research has focused largely on the seminorms of Gagliardo as an object capturing the fractional differential energy (see, for example, the Hitchhiker's Guide to the Fractional Sobolev spaces [12] for a thorough treatment). We also mention the interesting work of Kuusi and Mingione [20], who use the fractional maximal function applied to the gradient as an heuristic for fractional differentiation in the development of regularity estimates for *integer* order nonlinear elliptic equations, showing how a notion of fractional differentiation is useful even when one is interested in classical equations.

As there will be no ambiguity, in what follows we will contract distributional Riesz fractional gradient to fractional gradient. The remainder of the introduction will be dedicated to supporting our claim that the fractional gradient is an object of intrinsic interest, a task we will accomplish by presenting a unified theory of spaces of weakly fractionally differentiable functions and fractional partial differential equations. As with the theory of fractional Sobolev spaces detailed in [12], this theory continuously interpolates the existing integer order Sobolev space theory (and in fact is a theory associated to the complex interpolation of $L^p$ and $W^{1,p}$, in contrast to the real interpolation spaces obtained with the Gagliardo semi-norms). As such, there are objects in the integer setting that give rise to parallel structures in the interpolated theories. In particular, the reader should keep this dichotomy in mind when we discuss our notion of linear fractional partial differential equations with variable coefficients or the fractional p-Laplace's equation - while it is not clear the former has a direct analogy in the real interpolation setting, the latter has been defined as the operator associated to minimizing the Gagliardo semi-norm, e.g. [4, 11, 22]. In the setting of interpolation, the real and complex method complement each other to provide a more full understanding of the theory, and similarly, we hope our work concerning fractional partial differential equations in the complex interpolation spaces can provide a more complete understanding of the theory of fractional partial differential equations.

---

[3] The authors would like to thank Moritz Kassmann for making them aware of this reference.
[4] Our analysis is focused in the regime $N \geq 2$, though analogous theorems can be deduced when $N = 1$.

One of the first elements of our theory is the following theorem, which provides us with a geometric interpretation of the fractional gradient.

**Theorem 1.2.** *Let* $s \in (0,1)$. *If* $u \in C_c^\infty(\mathbb{R}^N)$, *then*

$$D^s u = I_{1-s} * Du.$$

Whereas the gradient provides one with the locally linear map approximating a function, Theorem 1.2 informs us that the fractional gradient is an averaged linear approximation. In this sense, nonlocal information about the derivative is found embedded into the fractional gradient.

We have mentioned the prominence of study of problems concerning the fractional Laplacian. Our next theorem connects the fractional partial derivatives and fractional Laplacian in an intuitive way, transforming the geometric meaning of the fractional gradient into a geometric understanding of the fractional Laplacian.

**Theorem 1.3.** *Let* $s \in (0,1)$. *If* $u \in C_c^\infty(\mathbb{R}^N)$, *then*

$$(-\Delta)^s u = -\sum_{j=1}^N \frac{\partial^s}{\partial x_j^s} \frac{\partial^s}{\partial x_j^s} u.$$

Compared with the standard definition of the Laplacian on a manifold or in Euclidean space, the preceding theorem is quite natural.

Along with the preceding real space geometric understanding of the fractional gradient, we have the following theorem computing its Fourier transform, providing us with a spectral understanding.

**Theorem 1.4.** *Let $s \in (0,1)$. Then*

$$\widehat{\frac{\partial^s u}{\partial x_j^s}} = -(2\pi)^s i\xi_j |\xi|^{-1+s} \widehat{u} \qquad (7)$$

*for all $u \in \mathscr{S}(\mathbb{R}^N)$.*

The preceding theorems provide us with insight into the object of the fractional gradient and its connection with the fractional Laplacian, and so we now proceed to develop spaces of weakly fractionally differentiable functions that will serve as a foundation for our existence results in the calculus of variations and partial differential equations. We therefore define the following space of fractionally differentiable functions.

**Definition 1.5.** Let $1 < p < \infty$ and $s \in (0,1)$. If $u \in C_c^\infty(\mathbb{R}^N)$, we define the norm

$$\|u\|_{X^{s,p}(\mathbb{R}^N)}^p := \|u\|_{L^p(\mathbb{R}^N)}^p + \|D^s u\|_{L^p(\mathbb{R}^N)}^p,$$

and the space

$$X^{s,p}(\mathbb{R}^N) := \overline{\{u \in C_c^\infty(\mathbb{R}^N)\}}^{\|\cdot\|_{X^{s,p}(\mathbb{R}^N)}}.$$

**Remark 1.6.** We have restricted our attention to the case $s \in (0,1)$, though our analysis can be extended to the case $s \geq 0$ by considering $\tilde{s} = s - \lfloor s \rfloor$, where $\lfloor s \rfloor$ is the integer part of $s$, and defining the norm on $X^{s,p}(\mathbb{R}^N)$ to be

$$\|u\|_{X^{s,p}(\mathbb{R}^N)}^p := \|u\|_{W^{\lfloor s \rfloor,p}(\mathbb{R}^N)}^p + \|D^{\tilde{s}} D^{\lfloor s \rfloor} u\|_{L^p(\mathbb{R}^N)}.$$

Having defined a new functional space, it would be appropriate to prove a number of theorems concerning its properties. We accomplish this task with a single theorem, identifying the spaces $X^{s,p}(\mathbb{R}^N)$ with the Bessel potential spaces $L^{s,p}(\mathbb{R}^N)$ (see Section 2 for a definition).

**Theorem 1.7.** *If $1 < p < \infty$ and $s \in (0,1)$, then*

$$X^{s,p}(\mathbb{R}^N) = L^{s,p}(\mathbb{R}^N).$$

This result is comparable in the integer order case to a theorem of Calderón for $s = k \in \mathbb{N}_0$ and $1 < p < \infty$ (where $L^{k,p}(\mathbb{R}^N)$ coincides with the Sobolev spaces $W^{k,p}(\mathbb{R}^N)$), and in the fractional case to a theorem of Strichartz concerning the nonlinear fractional operators $\mathcal{S}^s$ (see Adams and Hedberg, Chapter 3, Section 5 [2], or Stein p. 161-163 [31]). An important consequence of Theorem 1.7 is that we can deduce many properties about $X^{s,p}(\mathbb{R}^N)$ from the existing literature, as we will record in Section 2. However, with the addition of the fractional gradient, it raises the question of further analogy to the integer order case. It turns out that our viewpoint allows us to write and prove Sobolev and Hardy inequalities with great ease. For example, when $1 < p < \frac{N}{s}$, we have the fractional Sobolev inequality.

**Theorem 1.8** (Fractional Sobolev Inequality). *Let $1 < p < \infty$ and $s \in (0,1)$ be such that $sp < N$. Then there exists $C = C(N, p, s) > 0$ such that*

$$\|u\|_{L^{p^*}(\mathbb{R}^N)} \leq C \|D^s u\|_{L^p(\mathbb{R}^N)}$$

*for all $u \in L^{s,p}(\mathbb{R}^N)$, where $p^* := \frac{Np}{N-sp}$.*

In this regime, we also have the fractional Hardy inequality.

**Theorem 1.9** (Fractional Hardy Inequality). *Let $1 < p < \infty$, $s \in (0,1)$ and suppose $sp < N$. Then there exists $C > 0$ such that for every $u \in L^{s,p}(\mathbb{R}^N)$ we have*

$$\int_{\mathbb{R}^N} \frac{|u|^p}{|x|^{sp}} \, dx \leq C \int_{\mathbb{R}^N} |D^s u|^p \, dx.$$

Our statement of the preceding theorems serves to demonstrate the complete analogy of the fractional gradient, which extends past the regime $sp < N$. When $sp = N$, we are in the critical exponent case, and it is known that $L^{s,p}(\mathbb{R}^N)$ embeds into $L^q_{loc}(\mathbb{R}^N)$ for all $q \geq 1$. Actually, Trudinger's embedding from the integer order case can be extended to show that we have an embedding into $BMO(\mathbb{R}^N)$.

**Theorem 1.10** (Fractional Trudinger Inequality). *Let $1 < p < \infty$ and $s \in (0,1)$ be such that $sp = N$. Then there exists $A_1, A_2, C > 0$ such that for all $\Omega \subset \mathbb{R}^N$ open with finite measure*

$$\fint_\Omega exp \left[ \frac{|u(x)|}{A_1 C \|D^s u\|_{L^p(\mathbb{R}^N)}} \right]^{p'} dx \leq A_2$$

*for all $u \in L^{s,p}(\mathbb{R}^N)$.*

Finally, when $sp > N$, we have the analogy to the classical embedding theorem of Morrey (see Theorem 8 in [12] for an analagous theorem where $\|D^s u\|_{L^p(\mathbb{R}^N)}$ is replaced by the Gagliardo semi-norms).

**Theorem 1.11** (Fractional Morrey Inequality). *Let $1 < p < \infty$ and $s \in (0,1)$ be such that $sp > N$. Then there exists $M = M(N, p, s) > 0$ such that*

$$|u(x) - u(y)| \leq M|x-y|^{s-\frac{N}{p}} \|D^s u\|_{L^p(\mathbb{R}^N)}$$

*for all $u \in L^{s,p}(\mathbb{R}^N)$.*

An interesting tool we use in proving these theorems is a formula that deserves recording in its own right, the fractional fundamental theorem of calculus. The case $s = 1$ can be found in Chapter 5, Section 2, Equation (17) of Stein [31].

**Theorem 1.12** (Fractional FTOC). *Let $1 \leq p < \infty$ and $s \in (0,1)$. Every $u \in C_c^\infty(\mathbb{R}^N)$ is expressible as*

$$u = I_s \sum_{j=1}^{N} \mathcal{R}_j \frac{\partial^s u}{\partial x_j^s},$$

*where $\mathcal{R}_j$ is the Riesz transform, which can be characterized as a singular integral or $0$-order operator with symbol $\frac{-i\xi_j}{|\xi|}$.*

We now state our main results about the existence of minimizers for integral functionals of the fractional gradient and of weak solutions to fractional partial differential equations. We will treat the linear setting first, where we have the following theorem establishing the existence of solutions to linear fractional partial differential equations with variable coefficients.

**Theorem 1.13.** *Let $\Omega \subset \mathbb{R}^N$ be open and bounded. Suppose that $g \in L^2(\Omega)$, $\varphi \in H^s(\mathbb{R}^N)$ such that $I_{1-s} * \varphi$ is well-defined, and $A : \mathbb{R}^N \to \mathbb{R}^{N \times N}$ is bounded and measurable such that*

$$\lambda |z|^2 \leq A(x)z \cdot z \leq \Lambda |z|^2$$

*for some $\lambda, \Lambda > 0$ and all $x \in \mathbb{R}^N$ and $z \in \mathbb{R}^N$. Then there exists a unique $u \in H_\varphi^s(\Omega)$ such that*

$$\int_{\mathbb{R}^N} A(x) D^s u \cdot D^s v \, dx = \int_\Omega gv \, dx \tag{8}$$

*for every $v \in H_0^s(\Omega)$.*

**Remark 1.14.** Here, $H^s(\mathbb{R}^N) = L^{s,2}(\mathbb{R}^N)$ (although it has a number of equivalent definitions), and $u \in H_\varphi^s(\Omega)$ means $u \in H^s(\mathbb{R}^N)$ and $u = \varphi$ in $\Omega^c$.

As a corollary, we obtain existence of solutions to the fractional Laplace's equation.

**Corollary 1.15.** *Let $\Omega \subset \mathbb{R}^N$ be open. Suppose that $g \in L^2(\Omega)$ and $\varphi \in H^s(\mathbb{R}^N)$. There exists a unique $u \in H^s_\varphi(\Omega)$ that satisfies*

$$\int_{\mathbb{R}^N} (-\Delta)^{\frac{s}{2}} u (-\Delta)^{\frac{s}{2}} v \, dx = \int_\Omega gv \, dx$$

*for every $v \in H^s_0(\Omega)$.*

In the more general nonlinear setting, we need to define appropriate energy spaces analogous to $H^s_\varphi(\Omega)$.

**Definition 1.16.** Let $\Omega \subset \mathbb{R}^N$ be open, $1 < p < \infty$, and $s \in (0,1)$, and $\varphi \in L^{s,p}(\mathbb{R}^N)$. We define

$$L^{s,p}_\varphi(\Omega) := \{u \in L^{s,p}(\mathbb{R}^N) : u = \varphi \text{ in } \Omega^c\}.$$

With the energy space defined, we can now state our theorem concerning the existence of minimizers of integral functionals of the fractional gradient. The hypothesis on the integrand $f$ are comparable with those typically assumed in standard existence theorems in the calculus of variations, see the references [10, 14].

**Theorem 1.17.** *Let $1 < p < \infty$, $s \in (0,1)$, $g \in L^q(\Omega)$ with $q \geq (p^*)'$ if $sp < N$, and $\varphi \in L^{s,p}(\mathbb{R}^N)$ such that $I_{1-s} * \varphi$ is well-defined. Here, $p^* := \frac{Np}{N-sp}$ and $(p^*)'$ is the dual exponent. Suppose $f : \mathbb{R}^N \times \mathbb{R}^N \to \mathbb{R}$ is lower semicontinuous and convex in the second variable, and that there exists $c, C > 0$ and $\gamma_1, \gamma_2 \in L^1(\mathbb{R}^N)$ such that*

$$c|z|^p + \gamma_1(x) \leq f(x,z) \leq \gamma_2(x) + C|z|^p$$

*for all $(x,z) \in \mathbb{R}^N \times \mathbb{R}^N$. Then the energy*

$$F_s(u) := \int_{\mathbb{R}^N} f(x, D^s u) \, dx - \int_\Omega gu \, dx$$

*attains its infimum over $L^{s,p}_\varphi(\Omega)$.*

Assuming further differentiability of $f$ and upper bounds on the growth of its gradient in the second variable, we obtain the following theorem concerning the existence of weak solutions to fractional partial differential equations.

**Theorem 1.18.** *Let $1 < p < \infty$, $s \in (0,1)$, and that $\varphi \in L^{s,p}(\mathbb{R}^N)$ such that $I_{1-s} * \varphi$ is well-defined. Suppose $f \in C^1(\mathbb{R}^N \times \mathbb{R}^N)$ satisfies the hypothesis of Theorem* 1.17, *and in addition that there exists $a \in L^{p'}(\mathbb{R}^N)$ and $C > 0$ such that*

$$|\nabla_z f(x,z)| \leq a(x) + C|z|^{p-1}.$$

*Then there exists $u \in L^{s,p}_\varphi(\Omega)$ such that*

$$\int_{\mathbb{R}^N} \nabla_z f(x, D^s u) \cdot D^s v = \int_\Omega gv\, dx$$

*for all $v \in C_c^\infty(\Omega)$.*

In this setting, we obtain solutions to a fractional p-Laplace's equation.

**Corollary 1.19.** *Let $1 < p < \infty$, $s \in (0,1)$, and that $\varphi \in L^{s,p}(\mathbb{R}^N)$ such that $I_{1-s} * \varphi$ is well-defined. There exists $u \in L^{s,p}_\varphi(\Omega)$ such that*

$$\int_{\mathbb{R}^N} |D^s u|^{p-2} D^s u \cdot D^s v = \int_\Omega gv\, dx$$

*for all $v \in C_c^\infty(\Omega)$.*

**Remark 1.20.** As mentioned in the introduction, this is not to be confused with existing notions of the fractional p-Laplace's equation from the real interpolation theory, see for example, [5, 11, 22].

The organization of the paper will be as follows. In Section 2, we record some preliminary results regarding Riesz and Bessel potentials that play an essential role in our treatment of the fractional gradient. In Section 3, we prove the claims in the introduction concerning the fractional gradient, spaces of weakly fractionally differentiable functions, and fractional Sobolev inequalities. In Section 4, we treat the linear theory of fractional partial differential equations in several dimensions, while in Section 5, we will demonstrate some existence results for nonlinear the problem of existence of minimizers for integral functionals of the fractional gradient and corresponding existence results for weak fractional partial differential equations.

## 2  Preliminaries

Let us now define the space $L^{s,p}(\mathbb{R}^N)$, as well as state some of its properties. We first recall the Bessel potentials $g_s$, for $s \in \mathbb{R}_+$. The Bessel potentials $g_s$ are

defined by (see [24])

$$g_s(x) := \frac{1}{(4\pi)^{\frac{s}{2}}\Gamma(\frac{s}{2})} \int_0^\infty e^{\frac{-\pi|x|^2}{\delta}} e^{\frac{-\delta}{4\pi}} \delta^{\frac{s-N}{2}} \frac{d\delta}{\delta},$$

and can be shown to satisfy, for $s, t > 0$

$$\hat{g}_s(\xi) = (1 + 4\pi^2|\xi|^2)^{\frac{-s}{2}},$$
$$\|g_s\|_{L^1(\mathbb{R}^N)} = 1,$$
$$g_s * g_t = g_{s+t}$$

Then the Bessel Potential spaces $L^{s,p}(\mathbb{R}^N)$ are defined as follows.

**Definition 2.1.** Let $1 < p < +\infty$ and $s \in \mathbb{R}_+$. Define

$$L^{s,p}(\mathbb{R}^N) := g_s(L^p(\mathbb{R}^N)),$$

in the sense that every $u \in L^{s,p}(\mathbb{R}^N)$ can be written as

$$u = g_s * f,$$

for some $f \in L^p(\mathbb{R}^N)$.

As Theorem 1.7 records the equivalence of $X^{s,p}(\mathbb{R}^N)$ and $L^{s,p}(\mathbb{R}^N)$, the following theorem from Chapter 7 p.221 in the book of Adams [1] allows us to understand the properties of this space.

**Theorem 2.2** (Properties of $L^{s,p}(\mathbb{R}^N)$). (a) *If $s \geq 0$ and $1 \leq p < \infty$, then $C_c^\infty(\mathbb{R}^N)$ is dense in $L^{s,p}(\mathbb{R}^N)$.*

(b) *If $1 < p < \infty$ and $p' = \frac{p}{p-1}$, then $[L^{s,p}(\mathbb{R}^N)]' = L^{-s,p'}(\mathbb{R}^N)$.*

(c) *If $t < s$, then $L^{s,p}(\mathbb{R}^N) \hookrightarrow L^{t,p}(\mathbb{R}^N)$.*

(d) *If $t \leq s$ and $1 < p \leq q \leq \frac{Np}{N-(s-t)p}$, then $L^{s,p}(\mathbb{R}^N) \hookrightarrow L^{t,q}(\mathbb{R}^N)$.*

(e) *If $0 \leq \mu \leq s - \frac{N}{p} < 1$, then $L^{s,p}(\mathbb{R}^N) \hookrightarrow C^{0,\mu}(\mathbb{R}^N)$.*

(f) *If $s$ is a non-negative integer and $1 < p < \infty$, then $L^{s,p}(\mathbb{R}^N)$ coincides with the space $W^{s,p}(\mathbb{R}^N)$, the norms in the two spaces being equivalent. This conclusion holds for any real $s \in (0, 1)$ if $p = 2$.*

(g) *If $1 < p < \infty$ and $\epsilon > 0$, then for every $s$ we have*

$$L^{s+\epsilon,p}(\mathbb{R}^N) \hookrightarrow W^{s,p}(\mathbb{R}^N) \hookrightarrow L^{s-\epsilon,p}(\mathbb{R}^N).$$

Here, the notation $\hookrightarrow$ signifies a continuous embedding with an inequality of the norms.

**Remark 2.3.** Item f records the result of Calderón mentioned in the introduction in the integer setting, and more generally, in combination with Theorem 1.7 shows that,
$$X^{s,2}(\mathbb{R}^N) = L^{s,2}(\mathbb{R}^N) = W^{s,2}(\mathbb{R}^N).$$

Theorem 2.2 informs us to expected continuous embeddings, while in Section 3 we will show some Sobolev inequalities related to these embeddings. In order to accomplish this purpose, we need several lemmas recording properties of the Riesz potential as a map between various spaces, which can be found in the book of Mizuta [24][Chapter 4, p.153-157] and Gilbarg and Trudinger[15][Chapter 7, p. 152-157]. The first result concerns the regime $sp < N$. Let us remark that the Riesz potential being well-defined requires good integrability at infinity, and for this $f \in L^q(\mathbb{R}^N)$ for some $1 \leq q < N$ suffices. This is, for instance, the case for compactly supported functions $f \in L^p(\mathbb{R}^N)$ for some $1 < p < \infty$. The following theorem can be found in [24][Section 4.2, Theorem 2.1, p. 153]

**Lemma 2.4.** *Let $1 < p < \infty$ and $s \in (0,1)$ be such that $sp < N$. Then for all $f \in L^p(\mathbb{R}^N)$ such that $I_s * f$ is well-defined, there exists $C = C(N,p,s) > 0$ such that*
$$\|I_s f\|_{L^{p^*}(\mathbb{R}^N)} \leq C\|f\|_{L^p(\mathbb{R}^N)},$$
*where $p^* := \frac{Np}{N-sp}$.*

When $sp = N$, the Riesz potential maps $L^p(\mathbb{R}^N)$ into $BMO(\mathbb{R}^N)$, made precise by the following exponential estimate, which can be found, for instance, in [24][Section 4.2, Theorem 2.3, p. 156].

**Lemma 2.5.** *Let $1 < p < \infty$ and $s \in (0,1)$ be such that $sp = N$. Then there exists $A_1, A_2 > 0$ such that for every $G \subset \mathbb{R}^N$ open and bounded we have*
$$\fint_G \exp\left[\frac{|I_s f(x)|}{A_1\|f\|_{L^p(\mathbb{R}^N)}}\right]^{p'} dx \leq A_2$$
*for all $f \in L^p(\mathbb{R}^N)$ such that $I_s * f$ is well-defined.*

When $sp > N$, we have the following theorem on the Hölder continuity of the Riesz potential of a function in $L^p(\mathbb{R}^N)$, which can be found in [24][Section 4.2, Theorem 2.2, p. 155]

**Lemma 2.6.** *Let $1 < p < \infty$ and $s \in (0,1)$ be such that $sp > N$. Suppose $f \in L^p(\mathbb{R}^N)$ and that $I_s * f$ is well-defined. Then $I^s f \in C^{0,\alpha}_{loc}(\mathbb{R}^N)$ and*

$$|I_s f(x) - I_s f(y)| \leq M|x-y|^{s-\frac{N}{p}} \|f\|_{L^p(\mathbb{R}^N)}$$

The last lemma gives local estimates for $I_s$ for any value of the ratio $\frac{sp}{N}$.

**Lemma 2.7.** *Let $1 < p < \infty$ and $s \in (0,1)$, and define $\frac{1}{p^*} := \frac{1}{p} - \frac{s}{N}$. If $\frac{1}{q} > \frac{1}{p^*}$ and $q \geq 1$, then there exists $C = C(G, N, s, p)$ such that*

$$\left( \int_G |I_s f|^q \, dx \right)^{\frac{1}{q}} \leq C \|f\|_{L^p(\mathbb{R}^N)}$$

*for all $f \in L^p(\mathbb{R}^N)$ such that $I_s * f$ is well-defined.*

Finally, we conclude this section with a lemma that is a Hardy Inequality for Riesz potentials due to Stein and Weiss [32].

**Lemma 2.8.** *Let $1 < p < \infty$, $s \in (0,1)$ and $sp - N < \gamma < N(p-1)$. Define $\mu := \gamma - sp$. Then there exists a $C > 0$ such that for every $f \in L^p(\mathbb{R}^N)$ such that $I_s * f$ is well-defined we have*

$$\int_{\mathbb{R}^N} |x|^\mu |I_s f|^p \, dx \leq C \int_{\mathbb{R}^N} |x|^\gamma |f|^p \, dx.$$

## 3 Connections to Classical Notions

In this section, we prove the lemmas and theorems in the introduction concerning the fractional gradient, as well as the results concerning spaces of fractionally differentiable functions. We begin by proving Theorem 1.2.

*Proof of Theorem* 1.2. Let $u, v \in C_c^\infty(\mathbb{R}^N)$, define $K_v := supp\, v$ and let $R > 0$ such that for all $x \in K_v$, $supp\, u(x - \cdot) \subset B(0, R)$. Then

$$\langle \frac{\partial^s u}{\partial x_j^s}, v \rangle = (-1) \langle I_{1-s} u, \frac{\partial v}{\partial x_j} \rangle,$$

$$= -\int_{\mathbb{R}^N} (I_{1-s} * u) \frac{\partial v}{\partial x_j} \, dx$$

$$= -\lim_{n \to \infty} \int_{\mathbb{R}^N} (I_{1-s} * u) \frac{v(x + h_n e_j) - v(x)}{h_n} \, dx$$

$$= -\lim_{n \to \infty} \int_{K_v} \int_{B(0,R)} \frac{1}{h_n} \frac{u(x - h_n e_j - y) - u(x - y)}{|y|^{N-1+s}} \, dy \, v(x) \, dx.$$

Now, for $x \in K_v$ and $y \in B(0, R)$ we have

$$\left|\frac{1}{h_n}\frac{u(x-h_n e_j - y) - u(x-y)}{|y|^{N-1+s}}\right| \leq \frac{L}{|y|^{N-1+s}} \in L^1(B(0,R)),$$

which implies that for every $x \in K_v$

$$\left|\int_{B(0,R)} \frac{1}{h_n}\frac{u(x-h_n e_j - y) - u(x-y)}{|y|^{N-1+s}}\, dy\right| \leq C(L, R).$$

Therefore, the pointwise almost everywhere convergence

$$\frac{1}{h_n}\frac{u(x-h_n e_j - y) - u(x-y)}{|y|^{N-1+s}} \to -\frac{\frac{\partial u}{\partial x_j}(x-y)}{|y|^{N-1+s}}$$

and Lebesgue's dominated convergence theorem implies

$$\langle \frac{\partial^s u}{\partial x_j^s}, v\rangle = \int_{\mathbb{R}^N} (I_{1-s} * \frac{\partial u}{\partial x_j})v\, dx,$$

and the desired result follows by localizing in $v$. □

We will now prove Theorem 1.4, computing the Fourier transform of the fractional gradient.

*Proof of Theorem* 1.4. Suppose that $v \in C_c^\infty(\mathbb{R}^N)$. Then by the definition of the Fourier transform of a tempered distribution and the definition of the fractional partial derivatives we have

$$\langle \widehat{\frac{\partial^s u}{\partial x_j^s}}, v\rangle := \langle \frac{\partial}{\partial x_j}(I_{1-s} * u), \widehat{v}\rangle = -\langle (I_{1-s} * u), \frac{\partial \widehat{v}}{\partial x_j}\rangle. \tag{9}$$

However, since $I_{1-s} \in \mathscr{S}(\mathbb{R}^N)'$ and $u \in \mathscr{S}(\mathbb{R}^N)$, we have that the convolution is well-defined as a tempered distribution and thus the Fourier transform of the convolution is the product of the Fourier transforms. Therefore, since $\widehat{I_{1-s}} = (2\pi|\xi|)^{-1+s}$ and $\frac{\partial \widehat{v}}{\partial x_j} = (2\pi i \xi_j v)\widehat{\phantom{x}}$, we conclude that

$$\langle \widehat{\frac{\partial^s u}{\partial x_j^s}}, v\rangle = -\langle (2\pi)^s i \xi_j |\xi|^{-1+s}\widehat{u}, v\rangle,$$

which is the desired result. □

**Remark 3.1.** Theorem 1.4 gives a simple proof of Theorem 1.7 in the case $p = 2$.

An important point in the development of our theory is the connection of the fractional partial derivatives with the fractional Laplacian, stated in Theorem 1.7, and formally proven in the following computation. Lemma 1.4 implies that

$$\widehat{\frac{\partial^s}{\partial x_j^s}} = -(2\pi)^s i\xi_j |\xi|^{-1+1s},$$

and so

$$-\sum_{j=1}^N (2\pi)^{2s} i^2 \xi_j^2 |\xi|^{-2+2s} \widehat{u} = (2\pi|\xi|)^{2s} \widehat{u},$$

which agrees with definition (1). This result could be made precise through density, though we prefer the following use of the distributional calculus, since it reveals the connection with definition (5).

*Proof of Theorem* 1.3. Suppose $u \in C_c^\infty(\mathbb{R}^N)$. Then $D^s u = I_{1-s} * Du \in L^p(\mathbb{R}^N)$ for some $1 < p < \infty$ and $I_{1-s} * D^s u = I_{2-2s} * Du$ is well-defined, and therefore $D^s D^s u$ is well-defined. Thus, for $v \in C_c^\infty(\mathbb{R}^N)$ one has

$$\langle \frac{\partial^s}{\partial x_j^s} \frac{\partial^s u}{\partial x_j^s}, v \rangle := -\int_{\mathbb{R}^N} I_{1-s} * (I_{1-s} * \frac{\partial u}{\partial x_j}) \frac{\partial v}{\partial x_j} \, dx$$

$$= -\int_{\mathbb{R}^N} I_{2-2s} * \frac{\partial u}{\partial x_j} \frac{\partial v}{\partial x_j} \, dx$$

by the semi-group property (3). Then $u \in C_c^\infty(\mathbb{R}^N)$, along with the argument of Theorem 1.2 allows us to justify the interchanging of order of differentiation and integration, so that

$$-\int_{\mathbb{R}^N} I_{2-2s} * \frac{\partial u}{\partial x_j} \frac{\partial v}{\partial x_j} \, dx = -\int_{\mathbb{R}^N} \frac{\partial}{\partial x_j} I_{2-2s} * u \frac{\partial v}{\partial x_j} \, dx$$

$$= \int_{\mathbb{R}^N} \frac{\partial^2}{\partial x_j^2} I_{2-2s} * uv \, dx,$$

where we have integrated by parts. Therefore, using the formula (4), we conclude

$$\langle -\sum_{j=1}^N \frac{\partial^s}{\partial x_j^s} \frac{\partial^s u}{\partial x_j^s}, v \rangle = \int_{\mathbb{R}^N} -\Delta I_{2-2s} * uv \, dx$$

$$= \int_{\mathbb{R}^N} I_{-2s} uv \, dx,$$

which agrees with the definition (5). □

We can now proceed to prove Theorem 1.7, connecting the spaces $X^{s,p}(\mathbb{R}^N)$ with the Bessel potential spaces.

*Proof of Theorem* 1.7. We first prove the inclusion $L^{s,p}(\mathbb{R}^N) \subset X^{s,p}(\mathbb{R}^N)$.

Suppose $u \in L^{s,p}(\mathbb{R}^N)$, which means that $u = g_s * f$ for some $f \in L^p(\mathbb{R}^N)$. Then by Young's inequality for convolutions we have

$$\|u\|_{L^p(\mathbb{R}^N)} \leq \|g_s\|_{L^1(\mathbb{R}^N)} \|f\|_{L^p(\mathbb{R}^N)},$$

and thus $u \in L^p(\mathbb{R}^N)$. It therefore remains to show that $D^s u \in L^p(\mathbb{R}^N)$ and the inclusion will be demonstrated. We begin by assuming $f \in C_c^\infty(\mathbb{R}^N)$, and then the general case will follow from density. Note that $f \in C_c^\infty(\mathbb{R}^N)$ and $u \in L^p(\mathbb{R}^N)$ implies that it makes sense to write

$$\frac{\partial^s u}{\partial x_j^s} = g_s * \frac{\partial^s f}{\partial x_j^s}.$$

Then computing the Fourier transform we have

$$\widehat{\frac{\partial^s u}{\partial x_j^s}} = (1 + 4\pi^2 |\xi|^2)^{\frac{-s}{2}} ((2\pi)^s i \xi_j |\xi|^{-1+s}) \widehat{f},$$

so that

$$\widehat{\frac{\partial^s u}{\partial x_j^s}} = \frac{-i\xi_j}{|\xi|} \frac{(2\pi|\xi|)^s}{(1 + 4\pi^2 |\xi|^2)^{\frac{s}{2}}} \widehat{f}.$$

But Chapter 7, Lemma 2.1 in [24] or Chapter V, Section 3.2, Lemma 2 in Stein [31], states that there exists a finite signed measure $\mu_s$ for which

$$\widehat{\mu_s}(\xi) = \frac{(2\pi|\xi|)^s}{(1 + 4\pi^2 |\xi|^2)^{\frac{s}{2}}},$$

and therefore inverting the Fourier transform we have

$$\frac{\partial^s u}{\partial x_j^s} = \mu_s * \mathcal{R}_j f,$$

where $\mathcal{R}_j$ is the Riesz transform. Another application of Young's inequality yields

$$\|D^s u\|_{L^p(\mathbb{R}^N)} \leq \|\mu_s\|_{M_b(\mathbb{R}^N)} \|\mathcal{R} f\|_{L^p(\mathbb{R}^N)},$$

where

$$\|\mu_s\|_{M_b(\mathbb{R}^N)} := \sup_{\phi \in C_c(\mathbb{R}^N), \|\phi\|_\infty \leq 1} \int_{\mathbb{R}^N} \phi \, d\mu_s.$$

The previous inequality, along with the boundedness of the Riesz transform $\mathcal{R} : L^p(\mathbb{R}^N) \to L^p(\mathbb{R}^N)$ for $1 < p < \infty$ implies the desired result.

Conversely, suppose $u \in X^{s,p}(\mathbb{R}^N)$. Define

$$f := \lambda_s * \left( u + \sum_{j=1}^N \mathcal{R}_j \frac{\partial^s u}{\partial x_j^s} \right),$$

where $\lambda_s$ is the finite signed measure for which

$$(1 + 4\pi^2 |\xi|^2)^{\frac{s}{2}} = \widehat{\lambda}_s [1 + (2\pi|\xi|)^s],$$

as proven in Lemma 2.2, [24]. Again applying Young's inequality and invoking boundedness of the Riesz transform on $L^p(\mathbb{R}^N)$ for $1 < p < \infty$ we obtain that $f \in L^p(\mathbb{R}^N)$. It remains to show that $u = g_s * f$. Supposing that $u \in \mathscr{S}(\mathbb{R}^N) \cap X^{s,p}(\mathbb{R}^N)$, we have that

$$\widehat{g_s f} = \widehat{g_s} \widehat{\lambda}_s \left( \widehat{u} + \sum_{j=1}^N \frac{-i\xi_j}{|\xi|} ((2\pi)^s i \xi_j |\xi|^{-1+s}) \widehat{u} \right).$$

Making simplifications in the summation on the right hand side, we have

$$\widehat{g_s f} = \widehat{g_s} \widehat{\lambda}_s \left( \widehat{u} + (2\pi|\xi|)^s \widehat{u} \right),$$

which by the definition of $\lambda_s$ implies

$$\widehat{g_s f} = \widehat{u},$$

and therefore $u = g_s * f$ whenever $u \in \mathscr{S}(\mathbb{R}^N)$. The general case then follows from density of $\mathscr{S}(\mathbb{R}^N)$ in $X^{s,p}(\mathbb{R}^N)$, since $C_c^\infty(\mathbb{R}^N) \subset \mathscr{S}(\mathbb{R}^N) \cap X^{s,p}(\mathbb{R}^N) \subset X^{s,p}(\mathbb{R}^N)$ and $X^{s,p}(\mathbb{R}^N)$ is defined as the completion of $C_c^\infty(\mathbb{R}^N)$ in the $X^{s,p}$ norm. □

An interesting fact to note is that we have here shown the spaces of weakly fractionally differentiable functions to coincide with the Bessel potential spaces, which are also equivalent to the spaces defined by complex interpolation of the Sobolev spaces. This is in contrast to the work on fractional Sobolev spaces defined via the Gagliardo seminorm, which coincide with the real interpolation of $L^p$ and $W^{1,p}$.

Let us now prove the fractional fundamental theorem of calculus, Theorem 1.12, an essential tool in our proofs of the Sobolev inqualities.

*Proof of Theorem* 1.12. As Theorem 1.4 proves $\widehat{\frac{\partial^s u}{\partial x_j^s}} = (2\pi)^s i\xi_j |\xi|^{-1+s} \widehat{u}$ and $\widehat{\mathcal{R}_j} = \frac{-i\xi_j}{|\xi|}$, the result follows from the identity $\widehat{I_s} = (2\pi|\xi|)^{-s}$ in $\mathscr{S}(\mathbb{R}^N)'$. □

As in the case of integer order spaces, there are varying regimes for Sobolev inequalities depending on the ratio of the exponent, derivative order and dimension. In the regime $sp < N$, we have an analogy to the inequality of Sobolev, Gagliardo, and Nirenberg, given in Theorem 1.8 in the Introduction.

*Proof of Theorem* 1.8. It suffices to prove the result for $u \in C_c^\infty(\mathbb{R}^N)$, and we thereby obtain the result for general $u \in L^{s,p}(\mathbb{R}^N)$ by a density argument. But for $u \in C_c^\infty(\mathbb{R}^N)$, we may invoke Theorem 1.12 to write $u = I_s g$, where

$$g = \sum_{j=1}^N \mathcal{R}_j \frac{\partial^s u}{\partial x_j^s},$$

and note that $g \in L^p(\mathbb{R}^N)$ and $sp < N$. Thus, by Lemma 2.4

$$\|I_s g\|_{L^{p^*}(\mathbb{R}^N)} \leq C \|g\|_{L^p(\mathbb{R}^N)}.$$

Then the assumptions $1 < p < \infty$ and the boundedness of $\mathcal{R}_j : L^p(\mathbb{R}^N) \to L^p(\mathbb{R}^N)$ implies

$$\|g\|_{L^p(\mathbb{R}^N)} \leq C \|D^s u\|_{L^p(\mathbb{R}^N)}, \tag{10}$$

which with the previous inequality gives the desired result. □

**Remark 3.2.** A more technical Sobolev inequality involving the Lorentz spaces can also be deduced for the fractional gradient, see the work [28, 30] for these estimates in the context of the fractional Laplacian. We have chosen to restrict attention to the less technical setting of $L^p$ spaces for convenience of exposition.

An advantage to our definition is that the same argument which gives the Sobolev inequality for $sp < N$ extends to the other two regimes. For example, when $sp = N$ we find that $L^{s,p}(\mathbb{R})$ embeds into $BMO(\mathbb{R}^N)$, which can be precisely captured in the following inequality analogous to that of Trudinger.

*Proof of Theorem* 1.10. The proof is identical to Theorem 1.8, except that here we are in the regime $sp = N$ where we must invoke Lemma 2.5 which gives the estimate

$$\fint_\Omega exp\left[\frac{|I_s g(x)|}{A_1 \|g\|_{L^p(\mathbb{R}^N)}}\right]^{p'} dx \leq A_2,$$

which along with inequality (10) and monotonicity of the exponential implies the desired result. □

*Proof of Theorem* 1.11. Again, the proof is identical, requiring modification only through the estimate for $sp > N$ utilizing Lemma 2.6. □

For the purpose of demonstrating existence results in subsequent sections, we require the following local Sobolev inequality.

**Theorem 3.3.** *Let* $1 < p < \infty$ *and* $s \in (0,1)$, *and define* $\frac{1}{p^*} := \frac{1}{p} - \frac{s}{N}$. *If* $\frac{1}{q} > \frac{1}{p^*}$ *and* $q \geq 1$, *then there exists* $C = C(\Omega, N, s, p)$ *such that*

$$\left(\int_\Omega |u|^q \, dx\right)^{\frac{1}{q}} \leq C \|D^s u\|_{L^p(\mathbb{R}^N)}$$

*for all* $u \in L^{s,p}(\mathbb{R}^N)$.

We conclude this section by proving the Hardy Inequality stated in the introduction, which proceeds in an analogous manner to the preceding proofs.

*Proof of Theorem* 1.9. The inequality for $u \in C_c^\infty(\mathbb{R}^N)$ is a consequence of the representation given in Theorem 1.12 and Lemma 2.8, and the general case follows from density. □

## 4 Linear Setting

In the linear setting, the full power of the calculus of variations is not necessary to obtain existence of solutions to fractional partial differential equations.[5] In fact, Theorem 1.7 connects the fractional Sobolev spaces $X^{s,2}(\mathbb{R}^N)$ to the fractional Sobolev spaces $H^s(\mathbb{R}^N)$ for which known compactness results can be invoked. For instance, the Hitchhiker's Guide [12] records the following theorem.

**Lemma 4.1.** *Let* $0 < s < \frac{N}{2}$. *There exists a* $C = C(N,s) > 0$ *such that*

$$\|u\|^2_{L^{2^*}(\mathbb{R}^N)} \leq C \int_{\mathbb{R}^N} \int_{\mathbb{R}^N} \frac{|u(x) - u(y)|^2}{|x - y|^{N+2s}} \, dy dx$$

*for all* $u \in H^s(\mathbb{R}^N)$ *with compact support, where* $2^* := \frac{2N}{N-2s}$.

---

[5] The authors would like to thank Itai Shafrir for pointing out this simplification in the linear setting.

This result, along with the Lax-Milgram theorem can then be used to prove the following theorem, establishing existence of solutions to linear fractional partial differential equations.

*Proof of Theorem* 1.13. One only has to rewrite equation (8) for $\tilde{u} = u - \varphi$ to consider the equation

$$B[\tilde{u}, v] = \int_{\mathbb{R}^N} gv - A(x)D^s\varphi \cdot D^s v \; dx, \tag{11}$$

where the bilinear mapping $B : H_0^s(\Omega) \times H_0^s(\Omega) \to \mathbb{R}$ is defined by

$$B[\tilde{u}, v] = \int_{\mathbb{R}^N} A(x)D^s\tilde{u} \cdot D^s v \; dx.$$

Then Theorem 1.4, Lemma 4.1, and Hölder's inequality imply that

$$B[\tilde{u}, \tilde{u}] \geq \lambda \int_{\mathbb{R}^N} |D^s\tilde{u}|^2 \; dx \geq c||\tilde{u}||_{H^s}^2$$

which shows that $B$ is coercive on $H_0^s(\Omega)$. Meanwhile, by the assumptions on $A$ and the Cauchy-Schwarz inequality we have

$$B[\tilde{u}, v] \leq (2\pi)^{2s}\Lambda ||u||_{H^s(\mathbb{R}^N)} ||v||_{H^s(\mathbb{R}^N)},$$

which shows continuity of $B$. Similar estimates imply that the map

$$v \mapsto \int_{\mathbb{R}^N} gv - A(x)D^s\varphi \cdot D^s v \; dx$$

is a bounded linear functional on $H_0^s(\Omega)$, and so one may apply the Lax-Milgram theorem to obtain existence of $\tilde{u} \in H_0^s(\Omega)$ which satisfies (11) so that $\tilde{u} + \varphi \in H_\varphi^s(\Omega)$ and satisfies (8). Uniqueness then follows from linearity, since if $u_1, u_2$ are two solutions to (8), then the function $w = u_1 - u_2 \in H_0^s(\Omega)$ and satisfies

$$\int_{\mathbb{R}^N} A(x)D^s w \cdot D^s v \; dx = 0$$

for every $v \in H_0^s$. Thus, letting $v = w$ one obtains

$$\int_{\mathbb{R}^N} |D^s w|^2 \; dx = 0,$$

which applying Theorem 1.4 and using the fact that $w = 0$ on $\mathbb{R}^N \setminus \Omega$ implies $w \equiv 0$. □

Finally, we prove Corollary 1.15.

*Proof.* This follows from taking $A(x) = cI$ for an appropriate choice of $c$, and making two applications of Parseval's Theorem. The first takes the fractional derivatives of the functions into Fourier space, while the second brings them back into real space as fractional Laplacians of order $\frac{s}{2}$. □

## 5 Fractional Calculus of Variations and Partial Differential Equations

In the Hilbert space setting, we were able to use the Lax-Milgram theorem to deduce existence of weak solutions to linear fractional partial differential equations with variable coefficients. In the nonlinear setting, which includes energies like $|z|^p$ where $p \neq 2$, we utilize a standard technique in the calculus of variations - the direct method of Tonelli.

*Proof of Theorem* 1.17. We define

$$C_s := \inf_{L^{s,p}_\varphi(\Omega)} \int_{\mathbb{R}^N} f(x, D^s u) \, dx - \int_\Omega gu \, dx.$$

Theorem 3.3 implies $\varphi \in L^q(\Omega)$ with $1 \leq q < +\infty$, while if $sp < N$, $1 \leq q \leq p^*$. Thus, we have the upper bound

$$C_s \leq C \int_{\mathbb{R}^N} |D^s \varphi|^p \, dx + \|g\|_{L^{(q)'}(\Omega)} \|\varphi\|_{L^q(\Omega)},$$

and so we know that $C_s < +\infty$. Thus, we may find a minimizing sequence $\{u_n\} \subset L^{s,p}_\varphi(\Omega)$ such that

$$C_s = \lim_{n \to \infty} \int_{\mathbb{R}^N} f(x, D^s u_n) - gu_n \, dx.$$

We claim that, up to a subsequence which we will not relabel,

$$u_n - \varphi \rightharpoonup u - \varphi \text{ in } L^q(\mathbb{R}^N)$$
$$D^s u_n \rightharpoonup v \text{ in } L^p(\mathbb{R}^N),$$

where $q = p^*$ if $sp < N$ and for every $1 < q < +\infty$ otherwise. We further claim that

$$v = D^s u.$$

If we can establish these claims, then we may pass the limit in the preceding equation to conclude

$$\lim_{n\to\infty}\left(\int_{\mathbb{R}^N} f(x, D^s u_n)\,dx - \int_\Omega g u_n\,dx\right) \geq \liminf_{n\to\infty}\int_{\mathbb{R}^N} f(x, D^s u_n) - \limsup_{n\to\infty}\int_\Omega g u_n\,dx$$

$$\geq \int_{\mathbb{R}^N} f(x, D^s u) - \int_\Omega g u\,dx,$$

where we have used lower semicontinuity of the first term with respect to weak convergence and continuity with respect to weak convergence of the second, see for example Theorem 6.54 in [14].

However, for $n$ large, using coercivity of $f$ and Young's inequality with $\epsilon$ we have that

$$C_s + 1 \geq c\int_{\mathbb{R}^N} |D^s u_n|^p\,dx + \int_{\mathbb{R}^N} \gamma_1(x)\,dx - C_\epsilon \|g\|_{L^{q'}(\Omega)} - \epsilon\|u_n\|_{L^q(\Omega)}.$$

By choosing $\epsilon$ small enough and using the Sobolev inequality from Theorem 3.3 to control the term $\epsilon\|u_n\|_{L^q(\Omega)}$ and even to replace it with a term with a good sign, we have

$$\|u_n\|_{L^q(\Omega)} + \|D^s u_n\|_{L^p(\mathbb{R}^N)} \leq \tilde{C}_\epsilon.$$

This inequality, along with the fact that $u_n = \varphi$ on $\Omega^c$ and $\varphi \in L^q(\Omega)$ implies that $\{u_n - \varphi\}_n$ is bounded in $L^q(\mathbb{R}^N)$ (and equals zero outside $\Omega$). Thus, these bounds and reflexivity of the $L^p$ spaces for $p > 1$ implies the existence of a subsequence such that $u_n - \varphi \rightharpoonup u - \varphi$ and $D^s u_n \rightharpoonup v$ as claimed. It therefore remains to verify $v = D^s u$. The subtlety here is that we would like to write

$$\langle \frac{\partial^s}{\partial x_j^s} u_n, w\rangle = -\int_{\mathbb{R}^N} u_n I_{1-s} * \frac{\partial w}{\partial x_j}\,dx,$$

and let $u_n \rightharpoonup u$ weakly in $L^q(\mathbb{R}^N)$, but firstly we do not know that $u_n \in L^q(\mathbb{R}^N)$, since we only have a local estimate, and secondly we do not a priori know that $I_{1-s} * \frac{\partial w}{\partial x_j}$ is in a good enough space to make sense of the product. We therefore proceed as follows. From the definition of $D^s$ we have for $w \in C_c^\infty(\mathbb{R}^N)$

$$\langle \frac{\partial^s}{\partial x_j^s}(u_n - \varphi), w\rangle = -\int_{\mathbb{R}^N} I_{1-s} * (u_n - \varphi)\frac{\partial w}{\partial x_j}\,dx$$

$$= -\int_\Omega (u_n - \varphi) I_{1-s} * \frac{\partial w}{\partial x_j}\,dx,$$

which makes sense since $supp\,(u_n - \varphi)$ is compact. Then the convergence $u_n - \varphi \rightharpoonup u - \varphi$ in $L^q(\mathbb{R}^N)$ implies

$$-\lim_{n\to\infty} \int_\Omega (u_n - \varphi) I_{1-s} * \frac{\partial w}{\partial x_j}\, dx = -\int_{\mathbb{R}^N} (u - \varphi) I_{1-s} * \frac{\partial w}{\partial x_j}\, dx,$$

which allows us to move the Riesz potential again across the integral on the right hand side to obtain

$$\lim_{n\to\infty} \langle \frac{\partial^s}{\partial x_j^s}(u_n - \varphi), w\rangle = -\int_{\mathbb{R}^N} I_{1-s} * u \frac{\partial w}{\partial x_j}\, dx + \int_{\mathbb{R}^N} I_{1-s} * \varphi \frac{\partial w}{\partial x_j}\, dx.$$

However, by linearity of the fractional derivatives we have

$$\langle \frac{\partial^s}{\partial x_j^s}(u_n - \varphi), w\rangle = \langle \frac{\partial^s}{\partial x_j^s} u_n, w\rangle - \langle \frac{\partial^s}{\partial x_j^s} \varphi, w\rangle,$$

and since $\frac{\partial^s u_n}{\partial x_j^s} \in L^p(\mathbb{R}^N)$, we have

$$\langle \frac{\partial^s u_n}{\partial x_j^s}, w\rangle = \int_{\mathbb{R}^N} \frac{\partial^s u_n}{\partial x_j^s} w\, dx.$$

Then using the weak convergence $\frac{\partial^s u_n}{\partial x_j^s} \rightharpoonup v_j$, we have

$$\lim_{n\to\infty} \langle \frac{\partial^s}{\partial x_j^s}(u_n - \varphi), w\rangle = \int_{\mathbb{R}^N} v_j w\, dx - \langle \frac{\partial^s}{\partial x_j^s} \varphi, w\rangle,$$

which when combined with the previous computation, and cancelling the fractional partial derivatives of $\varphi$ on both sides, yields the desired result. □

*Proof of Theorem* 1.18. If we can verify Gâteaux differentiability of $F_s$, then the proof is completed, since defining

$$I(t) := F_s(u + tv),$$

where $u$ is the minimizer of $F_s$ over $L^{s,p}_\varphi(\Omega)$ obtained in Theorem 1.17 and $v \in C_c^\infty(\Omega)$, then $I$ is differentiable and

$$I(0) = \min\{I(t) : t \in \mathbb{R}\}.$$

Thus,

$$I'(0) = \frac{d}{dt} F_s(u + tv) = \langle F'_s(u), v\rangle.$$

24    T.-T. Shieh and D. SpectorIt therefore remains to verify Gâteaux differentiability of $F_s$. However, we have

$$\begin{aligned}
|\langle F'_s(u), v\rangle| &= \left|\int_{\mathbb{R}^N} \nabla_z f(x, D^s u) \cdot D^s v \, dx - \int_\Omega gv\right| \\
&\leq \int_{\mathbb{R}^N} (a(x) + |D^s u|^{p-1})|D^s v| \, dx + \int_\Omega |gv| \, dx \\
&\leq (\|a\|_{L^{p'}(\mathbb{R}^N)} + \|D^s u\|^{p-1}_{L^p(\mathbb{R}^N)})\|D^s v\|_{L^p(\mathbb{R}^N)} + \|g\|_{L^{q'}(\Omega)}\|v\|_{L^q(\Omega)}.
\end{aligned}$$

This shows that $F_s$ is Gâteaux differentiable and the proof is complete. □

**Acknowledgments.** The first author would like to thank Professor Ming-Chih Lai for his support in doing this project. The second author would like to thank Itai Shafrir, Lorena Aguirre, Rahul Garg, Georgios Psaradakis and Baptiste Devyver for their insightful discussions during the preparation of this paper, as well as the CMMSC at National Chiao Tung University for its support during the undertaking of this research.

**Author information**

Tien-Tsan Shieh, Department of Applied Mathematics, National Chiao Tung University, Taiwan.
E-mail: `tshieh@math.nctu.edu.tw`

Daniel E. Spector, Department of Mathematics, Technion - Israel Institute of Technology, Israel.
E-mail: `dspector@tx.technion.ac.il`